\documentclass[11pt,a4paper]{article}
\usepackage{authblk}
\usepackage{hyperref}
\usepackage{amssymb}
\usepackage{amsthm}
\usepackage{amsmath}
\usepackage{multirow}
\usepackage{float}
\usepackage{graphicx}
\usepackage{listings}
\newtheorem{theorem}{Theorem}[section]

\newtheorem{lemma}[theorem]{Lemma}
\newtheorem{proposition}[theorem]{Proposition}

\newtheorem{definition}[theorem]{Definition}
\newtheorem{example}[theorem]{Example}

\newtheorem{remark}[theorem]{Remark}

\title{ \Large\bfseries Algorithms for Simultaneous Block Triangularization and Block Diagonalization of Sets of Matrices }
\author[1]{Ahmad Y. Al-Dweik}
\author[2]{Ryad Ghanam}
\author[3]{Gerard Thompson}
\author[4]{M. T. Mustafa}
\affil[1]{ Department of Mathematics, Statistics and Physics, College of Arts and Sciences, Qatar University, Doha, 2713, State of Qatar; aydweik@qu.edu.qa}
\affil[2]{Department of Liberal Arts $\&$ Sciences, Virginia Commonwealth University in Qatar, Doha 8095, Qatar; raghanam@vcu.edu}
\affil[3]{Department of Mathematics, University of Toledo, Toledo, OH 43606, USA; gerard.thompson@utoledo.edu}
\affil[4]{ Department of Mathematics, Statistics and Physics, College of Arts and Sciences, Qatar University, Doha, 2713, State of Qatar; tahr.mustafa@qu.edu.qa}
\begin{document}
\maketitle
\begin{abstract}
In a recent paper, a new method was proposed to find the common invariant subspaces of a set of matrices. This paper invstigates the more general problem of putting a set of matrices into  block triangular or block-diagonal form simultaneously. Based on common invariant subspaces, two algorithms  for simultaneous block triangularization and block diagonalization of sets of matrices are presented. As an alternate approach for simultaneous block diagonalization of sets of matrices by an invertible matrix,  a new algorithm is developed  based on the generalized eigen vectors of a commuting matrix. Moreover, a new characterization for the simultaneous block diagonalization by an invertible matrix is provided. The algorithms are applied to concrete examples using the symbolic manipulation system Maple.
\end{abstract}
\bigskip


\bigskip



\section{Introduction}
A problem that occurs frequently in a variety of mathematical contexts, is to find the common invariant subspaces of a single, or set of matrices. In the case of a single endomorphism or matrix, it is relatively easy to find all the invariant subspaces by using Jordan normal form. Also, some theoretical results are given only for the invariant subspaces of two matrices. However, when there are more than two matices, the problem becomes much harder and unexpected invariant subspaces may occur and there is no systematic method is known. In a recent article \cite{Dweik2021}, we have provided a new algorithms to determine common invariant subspaces of a single or of a set of matrices systematically.

In the present article we consider a more general version of this problem, that is, providing two algorithms  for simultaneous block triangularization and block diagonalization of sets of matrices. One of the main steps in the first two proposed algorithms, consists of finding the common invariant subspaces of matrices using  the new method proposed in the recent article \cite{Dweik2021}.
It is worth mentioning that an efficient algorithm to explicitly compute a transfer matrix which realizes the simultaneous block diagonalization of \emph{unitary} matrices whose decomposition in irreducible blocks (common invariant subspaces) \emph{is known from elsewhere} is given in \cite{Bischer2021}. An application of simultaneous block-diagonalization of
\emph{normal} matrices in quantum theory is presented in \cite{PKJ2016}.

In this article we shall be concerned with finite dimensions only. Of course the fact that a single \emph{complex} matrix can always be put into triangular form follows readily from the Jordan normal form Theorem \cite{Dummit2004}. For a set of matrices, Jacobson in  \cite{Jacobson1953} introduces the notion of a \emph{composition series} for a collection of matrices. The idea of a composition series for a group is quite familiar. One shows that any two such composition series for a given group, assuming that one such exists at all, are isomorphic, the Jordan-H\"older Theorem \cite{Dummit2004}. Jacobson  in  \cite{Jacobson1953} characterized the simultaneous block triangularization of a set of matrices by the existence of
 a chain $\{0\}=V_0 \subset V_1 \subset ... \subset V_t=\mathbb{C}^n$ of invariant subspaces with dimension $dim(V_i/V_{i-1})=n_i$. Therefore, in the context of a collection of matrices  $\Omega=\{A_i \}_{ i=1}^N$, the idea is to locate a common invariant subspace $V$ of \emph{minimal} dimension $d$ of a set of matrices $\Omega$. Assume $V$ is generated by the (linearly independent) set
$\mathcal{B}_1=\{u_1,u_2, ..., u_d \}$,  and let $\mathcal{B} = \{u_1,u_2, ..., u_d , u_{d+1},u_{d+2}, ..., u_n\}$ be a basis of $\mathbb{C}^n$ containing $\mathcal{B}_1$. Upon setting $S= (u_1,u_2, ..., u_d , u_{d+1},u_{d+2}, ..., u_n)$, then $S^{-1}A_i S$ has the block triangular form:
\begin{equation*}
S^{-1}A_i S=
\left(
{\begin{array}{cc}
  B_{1,1}^i  & B_{1,2}^i   \\
   0   & B_{2,2}^i  \\
\end{array} }
 \right),
\end{equation*}
for $i=1,...,n$. Thereafter, one may define a quotient of the ambient vector space and each of the
matrices in the given collection will pass to this quotient.
As such one defines
\begin{equation*}
T_i= B_{2,2}^i =
\begin{pmatrix}
 \textbf{0}_{(n-d)\times d} & \textbf{I}_{n-d} \end{pmatrix}S^{-1}A_i S \begin{pmatrix} \textbf{0}_{d \times (n-d)} \\ \textbf{I}_{n-d}
\end{pmatrix}.
\end{equation*}
Then one may begin again the process of looking for a common invariant subspace of \emph{minimal} dimension of a set of matrices $\{T_i \}_{ i=1}^N$ and iterate the procedure. Since all spaces and matrices are of finite dimension, the procedure must terminate at some point. Again, any two such composition series will be isomorphic.
When the various quotients and submatrices are lifted back to the original vector space, one obtains precisely block-triangular form for the original set of matrices. It is important to find a composition series in the construction in order to make the set of matrices as ``block-triangular as possible".

Dubi \cite{Dubi2009} gave an algorithmic approach to simultaneous triangularization of a set of matrices based on the idea of Jacobson in  \cite{Jacobson1953}. In the  case of simultaneous triangularization, it  can be understood as the existence of
a chain $\{0\}=V_0 \subset V_1 \subset ... \subset V_t=\mathbb{C}^n$ of invariant subspaces with dimension $dim(V_i)=i$.
We generalize his study to cover simultaneous \emph{block} triangularization of a set of matrices. The generalized algorithm depends on the novel algorithm for constructing invariant subspaces of a set of matrices given in  the recent article \cite{Dweik2021}.

Specht \cite{Specht1955} (see also \cite{Shapiro1979}) proved that if the associative algebra $\mathcal{L}$ generated by a set of matrices $\Omega$  over $\mathbb{C}$ satisfies $\mathcal{L}=\mathcal{L}^{*}$, then $\Omega$ admits simultaneous block triangularization if and only if it admits simultaneous block diagonalization, in both cases via a unitary matrix. Following a result of Specht,  we prove that a set of matrices  $\Omega$ admits simultaneous block diagonalization if and only if the set $\Gamma=\Omega \cup \Omega^{*}$ admits simultaneous block triangularization. Finally, an algorithmic approach to simultaneous block diagonalization of a set of matrices  based on this fact is proposed.

The latter part of this paper presents an alternate approach for simultaneous block diagonalization of a set of  $n \times n$ matrices $\{A_s\}_{ s=1}^N$ by an invertible matrix that does not require finding the common invariant subspaces. Maehara et al \cite{Maehara2011} introduced algorithm for simultaneous block diagonalization of a set of matrices by a unitary  matrix based on the existence of a Hermitian commuting matrix. Here we extend their algorithm to simultaneous block diagonalization of a set of matrices by an invertible  matrix based on the existence of a commuting matrix which is not necessarily Hermitian. For example, consider the set of matrices  $\Omega=\{A_i \}_{i=1}^2$  where
\begin{equation}
A_1= \left( \begin {array}{ccc} 1&0&0\\ \noalign{\medskip}2&2&0
\\ \noalign{\medskip}1&1&1\end {array}
\right) ,
A_2=\left( \begin {array}{ccc} 0&0&0\\ \noalign{\medskip}2&1&0
\\ \noalign{\medskip}0&1&0\end {array}
\right) .
\end{equation}
The only Hermitian matrix commuting with the set $\Omega$ is the identity matrix. Therefore, we can not apply the proposed algorithm given in \cite{Maehara2011}. However, one can verify that the following non Hermitian matrix $C$ commute with all the matrices $\{A_i \}_{ i=1}^2$
\begin{equation}
C= \left( \begin {array}{ccc} 0&0&0\\ \noalign{\medskip}2&1&0
\\ \noalign{\medskip}0&1&0\end {array}
\right).
\end{equation}
The matrix $C$  has distinct eigenvalues  $\lambda_1=0, \lambda_2=1$ with algebraic multiplicities $n_1=2, n_2=1$ respectively. Moreover, the matrix $C$ is not diagonalizable. Therfore, we can not construct the eigenvalue decomposition for the matrix $C$. However, one can decompose the matrix $C$  by its generalized eigen vectors as follows:
\begin{equation}
S^{-1}C S= \left(\begin {array}{ccc} 0&1&0\\ \noalign{\medskip}0&0&0
\\ \noalign{\medskip}0&0&1\end {array}
\right)=\left(\begin {array}{cc} 0&1\\ 0&0\\ \end {array}\right) \oplus \left(1\right),
\end{equation}
where
\begin{equation}
S= \left(\begin {array}{ccc} 0&-\frac{1}{2}&0\\ \noalign{\medskip}0&1&1
\\ \noalign{\medskip}1&0&1\end {array} 
\right).
\end{equation}
Initially, it is noted that the matrices $\{A_i \}_{ i=1}^2$ can be decomposed into two diagonal blocks by the constructed invertible matrix $S$ where
\begin{equation}
\begin {array}{cc}
S^{-1}A_1 S= \left(\begin {array}{cc} 1&\frac{1}{2}\\ 0&1\\ \end {array}\right) \oplus 
                      \left(2\right),&
S^{-1}A_2 S= \left(\begin {array}{cc} 0&1\\ 0&0\\ \end {array}\right) \oplus
                     \left(1\right).
\end {array}
\end{equation}
Then a new algorithm is developed for  simultaneous block diagonalization by an invertible matrix based on the generalized eigen vectors of a commuting matrix. Moreover, a new characterization is presented by proving that the existence of a commuting matrix that possess at least two distinct eigenvalues is the necessary and sufficient condition to gurantee the simultaneous block diagonalization by an invertible matrix.

An outline of the paper is as follows. In Section 2 we review several definitions pertaining to block-triangular and  block-diagonal matrices and state several elementary consequences that follow from them.  In Section 3, following a result of Specht \cite{Specht1955} (see also \cite{Shapiro1979}), we provide conditions for putting a set of matrices into  \emph{block-diagonal} form simultaneously. Furthermore, we apply the theortical results to provide two algorithms that enable a collection of matrices to be put into \emph{block-triangular} form or \emph{block-diagonal} form simultaneously by a unitary matrix based on the existence of invariant subspaces.
In Section 4, a new characterizationsit is presented by proving that the existence of a commuting matrix that possess at least two distinct eigenvalues is the necessary and sufficient condition to gurantee the simultaneous block diagonalization by an invertible matrix. Furthermore, we apply the theortical results to provide an algorithm that enable a collection of matrices to be put into  \emph{block-diagonal} form simultaneously by an invertible matrix based on the existence of a commuting matrix. Section 3 and 4 also provide concrete examples using the symbolic manipulation system Maple.
\section{Preliminaries}
Let  $\Omega$ be a set of $n \times n$ matrices over an algebraically closed field $\mathcal{F}$, and
let $\mathcal{L}$ denote the algebra generated by $\Omega$ over $\mathcal{F}$. Similarly, let  $\Omega^{*}$ be the set of  the conjugate transpose of each matrix in $\Omega$ and $\mathcal{L}^{*}$ denote the algebra generated by $\Omega^{*}$ over $\mathcal{F}$.
\begin{definition}
An $n \times n$ matrix $A$ is  $BT(n_1, . . . , n_t)$ provided $A$ is block
upper triangular with $t$ square blocks on the diagonal, of sizes $n_1, . . . , n_t$, where $t \geq 2$ and $n_1+. . . +n_t=n$. That is, a block upper triangular matrix $A$ has the form
\begin{equation}
{\mathbf{A}} =
\left(
{\begin{array}{*{20}c}
   {{\mathbf{A}}_{1,1} } & {{\mathbf{A}}_{1,2} } &  \cdots  & {{\mathbf{A}}_{1,t} } \\
   0 & {{\mathbf{A}}_{2,2} } &  \cdots  & {{\mathbf{A}}_{2,t} }  \\
    \vdots  &  \vdots  &  \ddots  &  \vdots   \\
   0 & 0 &  \cdots  & {{\mathbf{A}}_{t,t} }  \\
\end{array} }
 \right)
\end{equation}
where ${\mathbf{A}}_{i,j}$ is a square matrix for all $i = 1, ..., t$ and  $j = i, ..., t$.
\end{definition}
\begin{definition}
A set of $n \times n$ matrices $\Omega$ is $BT(n_1, . . . , n_t)$ if all of the matrices in $\Omega$ are
 $BT(n_1, . . . , n_t)$.
\end{definition}
\begin{remark}
A set of $n \times n$ matrices $\Omega$ admits a simultaneous triangularization if it is $BT(n_1, . . . , n_t)$ with $n_i=1$ for $i=1,...,t$.
\end{remark}
\begin{remark}
A set of $n \times n$ matrices $\Omega$ is  $BT(n_1, . . . , n_t)$ if and only if the algebra $\mathcal{L}$ generated by $\Omega$  is  $BT(n_1, . . . , n_t)$.
\end{remark}
\begin{proposition}\rm \cite{Specht1955}(see also \cite{Shapiro1979})
Let $\Omega$ be a nonempty set of complex $n \times n$ matrices. Then there is a nonsingular matrix $S$ such that
$S \Omega S^{-1}$ is $BT(n_1, . . . , n_t)$ if and only if there is a unitary matrix $U$ such that
$U \Omega U^{*}$ is $BT(n_1, . . . , n_t)$.
\end{proposition}
\begin{theorem}\label{th1}\rm \cite[Chapter IV]{Jacobson1953}
Let $\Omega$ be a nonempty set of complex $n \times n$ matrices. Then there is  a unitary matrix $U$ such that
$U  \Omega  U^{*}$ is $BT(n_1, . . . , n_t)$ if and only if the set $\Omega$ has a chain $\{0\}=V_0 \subset V_1 \subset ... \subset V_t=\mathbb{C}^n$ of invariant subspaces with dimension $dim(V_i/V_{i-1})=n_i$.
\end{theorem}
\begin{definition}
An $n \times n$ matrix $A$ is  $BD(n_1, . . . , n_t)$ provided $A$ is block
diagonal with $t$ square blocks on the diagonal, of sizes $n_1, . . . , n_t$, where $t \geq 2$, $n_1+. . . +n_t=n$ and the blocks off the diagonal are the zero matrices. That is, a block diagonal matrix $A$ has the form
\begin{equation}
{\mathbf{A}} =
\left(
{\begin{array}{*{20}c}
   {{\mathbf{A}}_1 } & 0 &  \cdots  & 0  \\
   0 & {{\mathbf{A}}_2 } &  \cdots  & 0  \\
    \vdots  &  \vdots  &  \ddots  &  \vdots   \\
   0 & 0 &  \cdots  & {{\mathbf{A}}_t }  \\
\end{array} }
 \right)
\end{equation}
where ${\mathbf{A}}_k$ is a square matrix for all $k = 1, ..., t$. In other words, matrix ${\mathbf{A}}$ is the direct sum of
${\mathbf{A}}_1, ..., {\mathbf{A}}_t$. It can also be indicated as ${\mathbf{A}}_{\text{1}}  \oplus {\mathbf{A}}_{\text{2}}  \oplus ... \oplus {\mathbf{A}}_{\text{t}}$.
\end{definition}
\begin{definition}
A set of $n \times n$ matrices $\Omega$ is $BD(n_1, . . . , n_t)$ if all of the matrices in $\Omega$ are
 $BD(n_1, . . . , n_t)$.
\end{definition}
\begin{remark}
A set of $n \times n$ matrices $\Omega$ admits a simultaneous diagonalization if it is $BD(n_1, . . . , n_t)$ with $n_i=1$ for $i=1,...,t$.
\end{remark}
\begin{remark}
A set of $n \times n$ matrices $\Omega$ is  $BD(n_1, . . . , n_t)$ if and only if the algebra $\mathcal{L}$ generated by $\Omega$  is  $BD(n_1, . . . , n_t)$.
\end{remark}
\begin{proposition}\label{pr2}\rm \cite{Specht1955}(see also \cite{Shapiro1979})
Let $\Omega$ be a nonempty set of complex $n \times n$ matrices and let $\mathcal{L}$ be the algebra generated by $\Omega$  over $\mathbb{C}$.  Suppose  $\mathcal{L}=\mathcal{L}^{*}$. Then there is a nonsingular matrix $S$ such that
$S \mathcal{L} S^{-1}$ is $BT(n_1, . . . , n_t)$ if and only if there is  a unitary matrix $U$ such that
$U \mathcal{L} U^{*}$ is $BD(n_1, . . . , n_t)$.
\end{proposition}
\section{Algorithms for simultaneous block triangularization and block diagonalization of a set of matrices based on the invariant subspaces}
 Dubi \cite{Dubi2009} gave an algorithmic approach to simultaneous \emph{triangularization} of a set of  $n \times n$ matrices. In this section, we will generalize his study to cover simultaneous \emph{block triangularization} and  simultaneous \emph{block diagonalization} of a set of  $n \times n$ matrices. The generalized algorithms depend on the novel algorithm for constructing invariant subspaces of a set of matrices given in  the recent article \cite{Dweik2021} and Theorem 3.3.
\begin{lemma}\label{lm1}\rm
Let $\Omega$ be a nonempty set of complex $n \times n$ matrices,  $\Omega^{*}$ be the set of  the conjugate transpose of each matrix in $\Omega$ and $\mathcal{L}$ be the algebra generated by $\Gamma=\Omega \cup \Omega^{*}$. Then $\mathcal{L}=\mathcal{L}^{*}$.
\end{lemma}
\proof
Let $A$ be a matrix in  $\in \mathcal{L}$. Then $A=P(B_1,...,B_m)$ for some multivariate polynomial $P(x_1,...,x_m)$ and matrices $\{B_i\}_{i=1}^m\in \Gamma$. Therefore, $A^{*}=P^*(B_1,...,B_m)=Q(B_1^*,...,B_m^*)$ for some multivariate polynomial $Q(x_1,...,x_m)$ where the matrices $\{B_i^*\}_{i=1}^m\in \Gamma^*=\Gamma$. Hence the matrix $A^* \in \mathcal{L}$
\endproof
\begin{lemma}\label{lm2}\rm
Let $\Omega$ be a nonempty set of complex $n \times n$ matrices and  $\Omega^{*}$ be the set of  the conjugate transpose of each matrix in $\Omega$ and $\Gamma=\Omega \cup \Omega^{*}$. Then there is  a unitary matrix $U$ such that $U \Gamma U^{*}$ is $BD(n_1, . . . , n_t)$ if and only if there is  a unitary matrix $U$ such that $U \Omega U^{*}$ is $BD(n_1, . . . , n_t)$.
\end{lemma}
\proof
Assume that  there exists  a unitary matrix $U$ such that $U \Omega U^{*}$ is $BD(n_1, . . . , n_t)$.
Then $(U \Omega U^{*})^{*}=U \Omega^{*} U^{*}$ is $BD(n_1, . . . , n_t)$. Hence, $U \Gamma U^{*}$ is $BD(n_1, . . . , n_t)$.
\endproof
\begin{theorem}\rm
Let $\Omega$ be a nonempty set of complex $n \times n$ matrices,  $\Omega^{*}$ be the set of  the conjugate transpose of each matrix in $\Omega$ and $\Gamma=\Omega \cup \Omega^{*}$. Then there is  a unitary matrix $U$ such that $U  \Omega  U^{*}$ is $BD(n_1, . . . , n_t)$ if and only if  there is a unitary matrix $U$ such that $U \Gamma U^{*}$  is $BT(n_1, . . . , n_t)$.
\end{theorem}
\proof
Let $\mathcal{L}$ be the algebra generated by $\Gamma$. Then $\mathcal{L}=\mathcal{L}^{*}$ using lemma \ref{lm1}.
Now by applying proposition \ref{pr2} and lemma \ref{lm2},  the following statements are equivalent :

There is a unitary matrix $U$ such that $U \Gamma U^{*}$  is $BT(n_1, . . . , n_t)$. \\
$\iff$ There is a unitary matrix $U$ such that $U \mathcal{L} U^{*}$  is $BT(n_1, . . . , n_t)$. \\
$\iff$ There is  a unitary matrix $U$ such that $U \mathcal{L} U^{*}$ is $BD(n_1, . . . , n_t)$. \\
$\iff$ There is  a unitary matrix $U$ such that $U \Gamma U^{*}$ is $BD(n_1, . . . , n_t)$.\\
$\iff$ There is  a unitary matrix $U$ such that $U \Omega U^{*}$ is $BD(n_1, . . . , n_t)$.
\endproof
\subsection{Algorithm $A$:  Simultaneous block triangularization of a set of matrices $n \times n$ matrices $\{A_i \}_{ i=1}^N$.}
\begin{enumerate}
\item Input: the set $\Omega=\{A_i \}_{ i=1}^N$.
\item Set $k=0, \mathcal{B}=\phi, s=n, T_i=A_i, S_2=I$.
\item Search for  a $d$-dimensional invariant subspace  $V=\langle v_1, v_2, ..., v_d \rangle$ of a set of matrices $\{T_i \}_{ i=1}^N$ starting from $d=1$ up to $d=s-1$. If one does not exists and $k=0$, abort and print “no simultaneous block triangularization”. Else, if one does not exists and $k\ne 0$, go to step (8). Else,  go to next step.
\item Set $V_{k+1}=(S_2 v_1~S_2 v_2~ ...~S_2 v_d), \mathcal{B}= \mathcal{B} \cup \{S_2 v_1, S_2 v_2, ..., S_2 v_d\}, S_1=( V_1~V_2~ ...~ V_{k+1} )$.
\item Find a basis  $\{u_1,u_2, ..., u_l \}$ for the orthogonal complement of  $\mathcal{B}$.
\item  Set $S_2 =(u_1~u_2~ ...~u_l ), S=(S_1~S_2)$
and:\\
$T_i=\begin{pmatrix} \textbf{0}_{(s-d)\times d} & \textbf{I}_{s-d} \end{pmatrix}S^{-1}A_i S \begin{pmatrix} \textbf{0}_{d \times (s-d)} \\ \textbf{I}_{s-d} \end{pmatrix}$.
\item Set $k=k+1, s=s-d$ and return to step (3).
\item Compute the QR decomposition of the invertible  matrix $S$, by means of the Gram–Schmidt process, to convert it to a unitary matrix $Q$.
\item Output: a unitary matrix $U$ as the conjugate transpose of the resulted matrix $Q$.
\end{enumerate}
\begin{remark}
If one uses any non-orthogonal complement in step 5 of Algorithm $A$, then the matrix $S$ is invertible such that  $S^{-1} \Omega S$  is $BT(n_1, . . . , n_t)$. But in such a case, one can not guarantee that $U  \Omega  U^{*}$ is $BT(n_1, . . . , n_t)$.
\end{remark}
\newcounter{example}
\begin{example}
The set of matrices  $\Omega=\{A_i \}_{i=1}^2$ admits simultaneous block triangularization where
\begin{equation}
A_1= \left( \begin {array}{cccccc}
3&2&1&0&1&1\\ 
0&5&0&0&0&0\\ 
0&1&4&0&1&2\\ 
1&3&1&1&1&3\\ 
0&2&0&0&2&5\\
0&1&0&0&0&6
\end {array} \right) ,
A_2=\left( \begin {array}{cccccc}
44&12&4&-4&8&4\\ 
0&36&0&0&0&-1\\ 
0&12&32&0&4&4\\ 
4&16&8&52&4&4\\ 
0&4&-1&0&28&8\\ 
0&4&0&0&0
&40\end {array}\right) .
\end{equation}
Applying Algorithm $A$ to the set $\Omega$ can be summarizes as follows:
\begin{itemize}
\item Input: $\Omega$.
\item Initiation step:\\
We have $k=0, \mathcal{B}=\phi, s=6, T_1=A_1, T_2=A_2, S_2=I$.
\item In the first iteration:\\
We found two-dimensional invariant subspace  $V=\langle e_1, e_4 \rangle$ of a set of matrices $\{T_i \}_{ i=1}^2$.
Therefore, $\mathcal{B}=\{e_1, e_4\}, S_1=(e_1, e_4), S_2=(e_2,e_3,e_5,e_6)$,
\begin{equation}
T_1= \left( \begin {array}{cccc}
5&0&0&0\\ 
1&4&1&2\\ 
2&0&2&5\\
1&0&0&6
\end {array}\right) ,
T_2=\left(\begin {array}{cccc}
36&0&0&-1\\ 
12&32&4&4\\
4&-1&28&8\\
4&0&0&40
\end {array} \right),
\end{equation}
$k=1$  and $s=4$.
\item In the second iteration:
We found two-dimensional invariant subspace  $V=\langle e_2,e_3 \rangle$ of a set of matrices $\{T_i \}_{ i=1}^2$.
Therefore, $\mathcal{B}=\{e_1, e_4, e_3, e_5\}, S_1=(e_1, e_4, e_3, e_5), S_2=(e_2,e_6)$,
\begin{equation}
T_1= \left(  \begin {array}{cc}
5&0\\ 
1&6
\end {array} \right) ,
T_2=\left(\begin {array}{cc}
36&-1\\
4&40\end {array}\right),
\end{equation}
$k=2$  and $s=2$.
\item In the third iteration:
There is no one-dimensional invariant subspace of a set of matrices  $\{T_i \}_{ i=1}^2$.
Therefore, $S=(e_1~ e_4~e_3~e_5~e_2~e_6)$, and the corresponding unitary matrix is
$$U= \left(\begin {array}{cccccc}
1&0&0&0&0&0\\
0&0&0&1&0&0\\ 
0&0&1&0&0&0\\
0&0&0&0&1&0\\ 
0&1&0&0&0&0\\ 
0&0&0&0&0&1
\end {array} \right)
$$
such that the set $U \Omega  U^{*}=\{U  A_i  U^{*}\}_{i=1}^2$ is $BT(2,2,2)$ where
\begin{equation}
\begin{array}{l}
U  A_1  U^{*}=\left( \begin {array}{cc|cc|cc}
 3&0&1&1&2&1\\ 
1&1&1&1&3&3\\
\hline
0&0&4&1&1&2\\
0&0&0&2&2&5\\
\hline
0&0&0&0&5&0\\ 
0&0&0&0&1&6
\end {array}\right) ,\\
U  A_2  U^{*}=\left( \begin {array}{cc|cc|cc}
44&-4&4&8&12&4\\ 
4&52&8&4&16&4\\
\hline 
0&0&32&4&12&4\\ 
0&0&-1&28&4&8\\
\hline
0&0&0&0&36&-1\\
0&0&0&0&4&40
\end {array}\right) .\\
\end{array}
\end{equation}
\end{itemize}
\end{example}
\subsection{Algorithm $B$:  Simultaneous block diagonalization of a set of matrices $n \times n$ matrices $\{A_i \}_{ i=1}^N$.}
\begin{enumerate}
\item Input: the set $\Omega=\{A_i \}_{ i=1}^N$.
\item  Construct the set $\Gamma=\Omega \cup \Omega^{*}$.
\item  Find a unitary matrix $U$ such that $U  \Gamma  U^{*}$ is $BT(n_1, . . . , n_t)$ using Algorithm $A$  .
\item Output: a unitary matrix $U$.
\end{enumerate}
\begin{remark}
Algorithm $B$ provides the finest block-diagonalization. Moreover, the number of the blocks equals the number the of the  invariant subspaces and the size of each block is $n_i \times n_i$ where $n_i$ is  the dimension of the invariant subspace.
\end{remark}
\begin{example}
The set of matrices  $\Omega=\{A_i \}_{i=1}^2$ admits simultaneous block diagonalization where
\begin{equation}
A_1= \left( \begin {array}{ccccccc}
3&0&0&0&0&0&0\\ 
0&2&0&0&0&0&0\\ 
0&0&2&0&0&0&0\\ 
0&0&0&1&0&0&0\\ 
0&0&0&0&1&0&0\\ 
0&0&0&0&0&1&0\\ 
0&0&0&0&0&0&3\end {array} \right) ,
A_2=\left( \begin {array}{ccccccc}
 0&0&0&0&0&0&0\\ 
0&0&0&0&0&0&0\\ 
0&1&0&0&0&0&0\\ 
0&0&0&0&0&0&0\\ 
0&0&0&0&0&0&0\\ 
0&0&0&1&0&0&0\\ 
1&0&0&0&0&0&0\end {array} \right) .
\end{equation}
Applying Algorithm $B$ to the set $\Omega$ can be summarizes as follows:
\begin{itemize}
\item Input: $\Gamma=\Omega \cup \Omega^{*}$.
\item Initiation step:\\
We have $k=0, \mathcal{B}=\phi, s=7, T_1=A_1, T_2=A_2, T_3=A_2^T, S_2=I$.
\item In the first iteration:\\
We found one-dimensional invariant subspace  $V=\langle e_5 \rangle$ of a set of matrices $\{T_i \}_{ i=1}^3$.
Therefore, $\mathcal{B}=\{e_5\}, S_1=(e_5), S_2=(e_1,e_2,e_3,e_4,e_6,e_7)$,
\begin{equation}
T_1= \left( \begin {array}{cccccc}
3&0&0&0&0&0\\ 
0&2&0&0&0&0\\ 
0&0&2&0&0&0\\ 
0&0&0&1&0&0\\ 
0&0&0&0&1&0\\ 
0&0&0&0&0&3
\end {array} \right) ,
T_2=\left( \begin {array}{cccccc}
0&0&0&0&0&0\\ 
0&0&0&0&0&0\\ 
0&1&0&0&0&0\\ 
0&0&0&0&0&0\\ 
0&0&0&1&0&0\\ 
1&0&0&0&0&0
\end {array} \right),
T_3=T_2^T,
\end{equation}
$k=1$  and $s=6$.
\item In the second iteration:
We found two-dimensional invariant subspace  $V=\langle e_4,e_5 \rangle$ of a set of matrices $\{T_i \}_{ i=1}^3$.
Therefore, $\mathcal{B}=\{e_5, e_4, e_6\}, S_1=(e_5~e_4~e_6), S_2=(e_1,e_2,e_3,e_7)$,
\begin{equation}
T_1= \left( \begin {array}{cccc}
3&0&0&0\\ 
0&2&0&0\\ 
0&0&2&0\\ 
0&0&0&3
\end {array} \right) ,
T_2=\left( \begin {array}{cccc}
 0&0&0&0\\
0&0&0&0\\
0&1&0&0\\ 
1&0&0&0
\end {array} \right),
T_3=T_2^T,
\end{equation}
$k=2$  and $s=4$.
\item In the third iteration:
We found two-dimensional invariant subspace  $V=\langle e_2, e_3 \rangle$ of a set of matrices $\{T_i \}_{ i=1}^3$.
Therefore, $\mathcal{B}=\{e_5, e_4, e_6, e_2, e_3\}, S_1=(e_5~e_4~e_6~e_2~e_3), S_2=(e_1,e_7)$,
\begin{equation}
T_1= \left( \begin {array}{cc}
3&0\\ \noalign{\medskip}
0&3
\end {array} \right) ,
T_2=\left( \begin {array}{cc}
0&0\\ \noalign{\medskip}
1&0
\end {array} \right),
T_3=\left( \begin {array}{cc}
0&1\\ \noalign{\medskip}
0&0
\end {array} \right),
\end{equation}
$k=3$  and $s=2$.
\item In the fourth iteration:
There is no one-dimensional invariant subspace of a set of matrices  $\{T_i \}_{ i=1}^3$.
Therefore, $S=(e_5~e_4~e_6~e_2~e_3~e_1~e_7)$, and the corresponding unitary matrix is
$$U= \left( \begin {array}{ccccccc}
0&0&0&0&1&0&0\\ 
0&0&0&1&0&0&0\\ 
0&0&0&0&0&1&0\\ 
0&1&0&0&0&0&0\\ 
0&0&1&0&0&0&0\\ 
1&0&0&0&0&0&0\\ 
0&0&0&0&0&0&1
\end {array} \right)
$$
such that the set $U \Omega  U^{*}=\{U  A_i  U^{*}\}_{i=1}^2$ is $BD(1,2,2,2)$ where
\begin{equation}
\begin{array}{l}
U  A_1  U^{*}=\left( \begin {array}{c} 1\end {array} \right)\oplus
\left( \begin {array}{cc} 1&0\\ \noalign{\medskip}0&1\end {array} \right) \oplus
 \left( \begin {array}{cc} 2&0\\ \noalign{\medskip}0&2\end {array} \right) \oplus
\left( \begin {array}{cc} 3&0\\ \noalign{\medskip}0&3\end {array} \right) ,\\
U  A_2  U^{*}=\left( \begin {array}{c} 0\end {array} \right)\oplus
\left( \begin {array}{cc} 0&0\\ \noalign{\medskip}1&0\end {array} \right) \oplus
 \left( \begin {array}{cc} 0&0\\ \noalign{\medskip}1&0\end {array} \right) \oplus
\left( \begin {array}{cc} 0&0\\ \noalign{\medskip}1&0\end {array} \right) .\\
\end{array}
\end{equation}
\end{itemize}
\end{example}
\begin{example}
The set of matrices  $\Omega=\{A_i \}_{i=1}^2$ admits simultaneous block diagonalization where
\begin{equation}
A_1= \left( \begin {array}{ccccccc}
3&0&0&0&0&0&0\\ 
0&2&0&0&0&0&0\\ 
0&0&2&0&0&0&0\\
0&0&0&1&0&0&0\\ 
0&0&0&0&1&0&0\\ 
0&0&0&0&0&1&0\\ 
0&0&0&0&0&0&3
\end {array} \right) ,
A_2=\left( \begin {array}{ccccccc}
0&0&0&0&0&0&0\\ 
0&0&0&1&0&0&0\\ 
0&1&0&0&0&0&0\\ 
0&0&0&0&0&0&0\\ 
0&0&0&0&1&0&0\\ 
0&0&0&0&1&0&0\\ 
1&0&0&0&0&0&0\end {array} \right) .
\end{equation}
Similarly, applying Algorithm $B$ to the set $\Omega$ provides the matrix $S=(e_6~e_5~e_7~e_1~ e_3~e_2~e_4)$. Therefore, the corresponding unitary matrix is
$$U= \left( \begin {array}{ccccccc}
 0&0&0&0&0&1&0\\ 
0&0&0&0&1&0&0\\
0&0&0&0&0&0&1\\ 
1&0&0&0&0&0&0\\ 
0&0&1&0&0&0&0\\
0&1&0&0&0&0&0\\ 
0&0&0&1&0&0&0
\end {array} \right)
$$
such that the set $U \Omega  U^{*}=\{U  A_i  U^{*}\}_{i=1}^2$ is $BD(2,2,3)$ where
\begin{equation}
\begin{array}{l}
U  A_1  U^{*}= \left( \begin {array}{cc} 1&0\\ \noalign{\medskip}0&1\end {array}\right) \oplus
 \left( \begin {array}{cc} 3&0\\ \noalign{\medskip}0&3\end {array}\right)\oplus
 \left( \begin {array}{ccc} 2&0&0\\ \noalign{\medskip}0&2&0\\ \noalign{\medskip}0&0&1\end {array} \right) ,\\
U  A_2  U^{*}= \left( \begin {array}{cc} 0&1\\ \noalign{\medskip}0&1\end {array}\right) \oplus
 \left( \begin {array}{cc} 0&1\\ \noalign{\medskip}0&0\end {array}\right)\oplus
 \left( \begin {array}{ccc} 0&1&0\\ \noalign{\medskip}0&0&1\\ \noalign{\medskip}0&0&0\end {array} \right)  .\\
\end{array}
\end{equation}
\end{example}
\begin{example}
The set of matrices  $\Omega=\{A_i \}_{i=1}^3$ admits simultaneous block diagonalization where
\begin{equation}
\begin {array}{ll}
A_1= \left( \begin {array}{ccccccccc}
0&0&0&0&0&0&0&0&0\\ 
0&2&0&0&0&0&0&0&0\\ 
0&0&1&0&0&0&0&0&0\\ 
0&0&0&-2&0&0&0&0&0\\ 
0&0&0&0&0&0&0&0&0\\ 
0&0&0&0&0&-1&0&0&0\\ 
0&0&0&0&0&0&-1&0&0\\
0&0&0&0&0&0&0&1&0\\ 
0&0&0&0&0&0&0&0&0
\end {array} \right)  ,
A_2=\left( \begin {array}{ccccccccc}
 0&0&0&1&0&0&0&0&0\\ 
-1&0&0&0&1&0&0&0&0\\ 
0&0&0&0&0&1&0&0&0\\ 
0&0&0&0&0&0&0&0&0\\ 
0&0&0&-1&0&0&0&0&0\\ 
0&0&0&0&0&0&0&0&0\\ 
0&0&0&0&0&0&0&0&0\\ 
0&0&0&0&0&0&-1&0&0\\ 
0&0&0&0&0&0&0&0&0
\end {array} \right),\\
A_3=\left( \begin {array}{ccccccccc}
 0&-1&0&0&0&0&0&0&0\\ 
0&0&0&0&0&0&0&0&0\\ 
0&0&0&0&0&0&0&0&0\\ 
1&0&0&0&-1&0&0&0&0\\
0&1&0&0&0&0&0&0&0\\
0&0&1&0&0&0&0&0&0\\ 
0&0&0&0&0&0&0&-1&0\\
0&0&0&0&0&0&0&0&0\\ 
0&0&0&0&0&0&0&0&0
\end {array} \right).
\end{array}
\end{equation}
Similarly, applying Algorithm $B$ to the set $\Omega$ provides the matrix $S=(e_1+e_5~e_9~e_3~e_6~ e_8~-e_7~ e_1-e_5,e_2~e_4)$. Therefore, the corresponding unitary matrix is
$$U= \left( \begin {array}{ccccccccc}
 \frac{1}{2\sqrt{2}}&0&0&0& \frac{1}{2\sqrt{2}}&0&0&0&0\\ 
0&0&0&0&0&0&0&0&1\\ 
0&0&1&0&0&0&0&0&0\\ 
0&0&0&0&0&1&0&0&0\\ 
0&0&0&0&0&0&0&1&0\\ 
0&0&0&0&0&0&-1&0&0\\ 
 \frac{1}{2\sqrt{2}}&0&0&0&- \frac{1}{2\sqrt{2}}&0&0&0&0\\ 
0&1&0&0&0&0&0&0&0\\ 
0&0&0&1&0&0&0&0&0
\end {array} \right) $$
such that the set $U \Omega  U^{*}=\{U  A_i  U^{*}\}_{i=1}^3$ is $BD(1,1,2,2,3)$ where
\begin{equation}
\begin{array}{l}
U  A_1  U^{*}=\left( \begin {array}{c} 0\end {array} \right)\oplus \left( \begin {array}{c} 0\end {array} \right) \oplus
 \left( \begin {array}{cc} 1&0\\ \noalign{\medskip}0&-1\end {array}\right) \oplus
 \left( \begin {array}{cc} 1&0\\ \noalign{\medskip}0&-1\end {array}\right)\oplus
 \left( \begin {array}{ccc} 0&0&0\\ \noalign{\medskip}0&2&0\\ \noalign{\medskip}0&0&-2\end {array} \right) ,\\
U  A_2  U^{*}=\left( \begin {array}{c} 0\end {array} \right)\oplus \left( \begin {array}{c} 0\end {array} \right) \oplus
 \left( \begin {array}{cc} 0&1\\ \noalign{\medskip}0&0\end {array}\right) \oplus
 \left( \begin {array}{cc} 0&1\\ \noalign{\medskip}0&0\end {array}\right) \oplus
 \left( \begin {array}{ccc} 0&0&\sqrt{2}\\ \noalign{\medskip}-\sqrt{2}&0&0\\ \noalign{\medskip}0&0&0\end {array} \right) ,\\
U  A_3  U^{*}=\left( \begin {array}{c} 0\end {array} \right)\oplus \left( \begin {array}{c} 0\end {array} \right) \oplus
 \left( \begin {array}{cc} 0&0\\ \noalign{\medskip}1&0\end {array} \right) \oplus
 \left( \begin {array}{cc} 0&0\\ \noalign{\medskip}1&0\end {array} \right) \oplus
\left( \begin {array}{ccc} 0&-\sqrt{2}&0\\ \noalign{\medskip}0&0&0\\ \noalign{\medskip}\sqrt{2}&0&0\end {array} \right) .\\
\end{array}
\end{equation}
\end{example}
\section{Algorithm for simultaneous  block diagonalization of a set of matrices based on a  commuting matrix}
This section focuses on an alternate approach for simultaneous block diagonalization of a set of  $n \times n$ matrices $\{A_s\}_{ s=1}^N$ by an invertible matrix that does not require finding the common invariant subspaces as Algorithm $B$ given in the previous section.
Maehara et al \cite{Maehara2011} introduced algorithm for simultaneous block diagonalization of a set of matrices by a unitary  matrix based on the eigenvalue decomposition of a Hermitian commuting matrix. Here we extend their algorithm to be applicable for non-Hermitian commuting matrix by considering its generalized eigen vectors. Moreover, a new characterizationsit is presented by proving that the existence of a commuting matrix that possess at least two distinct eigenvalues is the necessary and sufficient condition to gurantee the simultaneous block diagonalization by an invertible matrix.
\begin{proposition}\label{Pro1}\rm
Let $V$ be a vector space, and let $T:V \rightarrow V$ be a linear operator. Let $\lambda_1,..., \lambda_k$ be distinct eigenvalues of $T$. Then each 
generalized eigenspace $G_{\lambda_i}(T)$ is $T$-invariant, and we have  the direct sum decomposition
$$V=G_{\lambda_1}(T)\oplus G_{\lambda_2}(T) \oplus ...\oplus G_{\lambda_k}(T).$$
\end{proposition}
\begin{lemma}\label{Lem1}\rm
Let $V$ be a vector space, and let $T:V \rightarrow V$, $L:V \rightarrow V$ be linear commuting operators. Let $\lambda_1,..., \lambda_k$ be distinct eigenvalues of $T$. Then each  generalized eigenspace $G_{\lambda_i}(T)$ is $L$-invariant.
\end{lemma}
\proof
Let $V$ be a vector space, and $\lambda_1,..., \lambda_k$ be distinct eigenvalues of $T$ with the minimal polynomial $\mu(x)=(x-\lambda_1)^{n_1}(x-\lambda_2)^{n_2}...(x-\lambda_k)^{n_k}$. Then we have  the direct sum decomposition
$V=G_{\lambda_1}(T)\oplus G_{\lambda_2}(T) \oplus ...\oplus G_{\lambda_k}(T)$.
 
For each $i=1,..,k$, let  $x \in  G_{\lambda_i}(T)$, then $(T-\lambda_i I)^{n_i}x=0$.  Then  
$(T-\lambda_i I)^{n_i}Lx=L(T-\lambda_i I)^{n_i}x=0$. Hence $L x \in  G_{\lambda_i}(T)$.
\endproof
\begin{theorem}\rm
Let  $\{A_s\}_{ s=1}^N$ be a set of  $n \times n$ matrices. Then  the set  $\{A_s\}_{ s=1}^N$  admits simultaneous block diagonalization by an invertible matrix $S$ if and only if  the set  $\{A_s\}_{ s=1}^N$ commutes with  a matrix $C$ that possess two distinct eigenvalues. 
\end{theorem}
\proof
\begin{enumerate}
\item[$\Rightarrow$] 
Assume that the set  $\{A_s\}_{ s=1}^N$  admits simultaneous block diagonalization by the an invertible matrix $S$ such that
$$S^{-1}  A_s S=B_{s,1} \oplus  B_{s,2} \oplus ... \oplus B_{s,k},$$
where the number of blocks $k\geq 2$ and the matrices $B_{s,1}, B_{s,2}, ... ,B_{s,k}$  have sizes $n_1 \times n_1, n_2 \times n_2, ...,n_k \times n_k$ respectively for all $s=1,..,N$.

Now, define  the matrix $C$ as 
$$ C =S (\lambda_1 I_{n_1 \times n_1} \oplus  \lambda_2 I_{n_2 \times n_2} \oplus ... \oplus \lambda_k I_{n_k \times n_k})  S^{-1},$$
where  $\lambda_1,\lambda_2,..., \lambda_k$ are any distinct numbers. 

Clearly, the matrix $C$ commutes with  the set  $\{A_s\}_{ s=1}^N$. Moreover, it has the distinct eigenvalues $\lambda_1,\lambda_2,..., \lambda_k$.
\item[$\Leftarrow$] 
Assume that the set  $\{A_s\}_{ s=1}^N$ commutes with  a  matrix $C$ that posses distinct eigenvalues  $\lambda_1,\lambda_2,..., \lambda_k$. 

Using Proposition \ref{Pro1}, one can use the generalized eigenspace $G_{\lambda_i}(C)$  of the matrix $C$ associated to these distinct  eigenvalues to decompose the matrix $C$ as a direct sum of $k$ matrices. This can be achieved by restricting the matrix $C$ on the invariant subspaces $G_{\lambda_i}(C)$  as follows:
$$S^{-1}{C}S={\big[ C \big]}_{G_{\lambda_1}(C)} \oplus  {\big[ C \big]}_{G_{\lambda_2}(C)} \oplus  ... \oplus {\big[ C \big]}_{G_{\lambda_k}(C)}$$
where  $$S=\big( G_{\lambda_1}(C), G_{\lambda_2}(C) ,...,G_{\lambda_k}(C) \big).$$

Using Lemma \ref{Lem1}, one can restrict each matrix $A_s$ on the invariant subspaces $G_{\lambda_i}(C)$ to decompose the matrix $A_s$ as a direct sum of $k$ matrices as follows:
$$S^{-1}{A_s}S={\big[ A_s \big]}_{G_{\lambda_1}(C)} \oplus  {\big[ A_s \big]}_{G_{\lambda_2}(C)} \oplus  ... \oplus {\big[ A_s \big]}_{G_{\lambda_k}(C)}.$$
\end{enumerate}
\endproof
\begin{remark}
For a given set  of $n \times n$ matrices $\{A_s\}_{ s=1}^N$, if the set  $\{A_s\}_{ s=1}^N$ commutes only with the matrices having only one eigenvalue, then it does not admit a simultaneous block diagonalization by an invertible matrix.
\end{remark}
\subsection*{Algorithm $C$:}
\begin{enumerate}
\item Input: the set $\Omega=\{A_s \}_{ s=1}^N$. 
\item Construct the the following matrix:
\begin{equation*}
X=
 \left( \begin {array}{c}
I \otimes A_1  -A_1^T \otimes I  \\ \noalign{\medskip}
I \otimes A_2  -A_2^T \otimes I  \\ \noalign{\medskip}
.\\
.\\
.\\
I \otimes A_N  -A_N^T \otimes I  \\ \noalign{\medskip}
\end {array} \right) .
\end{equation*}
\item Compute the null space of the matrix $X$ and reshape the obtained vectors as $n \times n$ matrices. These matrices commute with all the matrices $\{A_s \}_{ s=1}^N$.
\item Choose a matrix $C$  from  the obtianed matrices that possess two distinct eigenvalues. 
\item Find the distinct eigenvalues $\lambda_1,..., \lambda_k$ of the matrix $C$ and  the corresponding algebraic multiplicity $n_1, n_2, ..., n_k$.
\item Find each generalized eigenspace $G_{\lambda_i}(C)$  of the matrix $C$ associated to the eigenvalue $\lambda_i$ by computing the null space of $(C-\lambda_i I)^{n_i}$.
\item Construct the invertible matrix $S$ as $$S=\big( G_{\lambda_1}(C), G_{\lambda_2}(C) ,...,G_{\lambda_k}(C) \big).$$
\item  Verify that
$$S^{-1}  A_s S=B_{s,1} \oplus  B_{s,2} \oplus ... \oplus B_{s,k},$$
where the matrices $B_{s,1}, B_{s,2}, ... ,B_{s,k}$  have sizes $n_1 \times n_1, n_2 \times n_2, ...,n_k \times n_k$ respectively for all $s=1,..,N$.
\item Output: an invertible matrix $S$.
\end{enumerate}
\begin{remark}
Algorithm $C$ provides the finest block-diagonalization if one choose a matrix $C$  with maximum number of distinct eigenvalues. Moreover, the number of the blocks equals the number the of the distinct eigenvalues and the size of each block is $n_i \times n_i$ where $n_i$ is  the algebraic multiplicity of the eigenvalue $\lambda_i$.
\end{remark}
\begin{example}
Consider the set of matrices  $\Omega=\{A_i \}_{i=1}^6$ where 
\begin{footnotesize}
\begin{equation}
\begin {array}{l}
A_1=\left( \begin {array}{cccccc}
0&0&0&0&0&0\\
0&0&0&1&0&0\\
0&0&0&0&1&0\\
0&-1&0&0&0&0\\
0&0&-1&0&0&0\\
0&0&0&0&0&0
\end {array} \right) , 
A_2=\left( \begin {array}{cccccc}
&0&0&-1&0&0\\
0&0&0&0&0&0\\
0&0&0&0&0&1\\
1&0&0&0&0&0\\
0&0&0&0&0&0\\
0&0&-1&0&0&0
\end {array} \right), 
A_3=\left( \begin {array}{cccccc}
 0&0&0&0&-1&0\\
0&0&0&0&0&-1\\
0&0&0&0&0&0\\
0&0&0&0&0&0\\
1&0&0&0&0&0\\ 
0&1&0&0&0&0
\end {array} \right) , \\
A_4=\left( \begin {array}{cccccc}
0&1&0&0&0&0\\ 
-1&0&0&0&0&0\\
0&0&0&0&0&0\\ 
0&0&0&0&0&0\\ 
0&0&0&0&0&1\\
0&0&0&0&-1&0
\end {array} \right), 
A_5=\left( \begin {array}{cccccc}
 0&0&1&0&0&0\\ 
0&0&0&0&0&0\\
-1&0&0&0&0&0\\
0&0&0&0&0&-1\\
0&0&0&0&0&0\\
0&0&0&1&0&0
\end {array} \right), 
A_6=\left( \begin {array}{cccccc}
0&0&0&0&0&0\\
0&0&1&0&0&0\\ 
0&-1&0&0&0&0\\ 
0&0&0&0&1&0\\
0&0&0&-1&0&0\\ 
0&0&0&0&0&0
\end {array} \right).
\end{array}
\end{equation}
\end{footnotesize}
The set $\Omega$  admits  simultaneous block diagonalization by an invertible matrix.  An invertible matrix can be obtained by applying algorithm $C$  to the set $\Omega$ as summarizes below:
\begin{itemize}
\item A matrix $C$  that commute with all the matrices $\{A_i \}_{ i=1}^6$ can be obtained as 
\begin{equation}
C= \left( \begin {array}{cccccc} 
0&0&0&0&0&1\\
0&0&0&0&-1&0\\
0&0&0&1&0&0\\ 
0&0&1&0&0&0\\
0&-1&0&0&0&0\\
1&0&0&0&0&0
\end {array} \right) 
\end{equation}.
\item The distinct eigenvalues of the matrix $C$ are $\lambda_1=-1, \lambda_2=1$ with algebraic multiplicities $n_1=3, n_2=3$ respectively..  
\item The generalized eigenspace  of the matrix $C$ associated to the distinct eigenvalues are
\begin{equation}
\begin {array}{l}
G_{\lambda_1}(C)=\mathcal{N}(C-\lambda_1 I)^3=\langle e_6-e_1,e_2+e_5,e_4-e_3\rangle,\\
G_{\lambda_2}(C)=\mathcal{N}(C-\lambda_2 I)^3=\langle e_1+e_6, e_5-e_2,e_3+e_4 \rangle.\\
\end {array} 
\end{equation}
\item  The invertible matrix $S=\big( G_{\lambda_1}(C),G_{\lambda_2}(C) \big)$ is 
\begin{equation}
S= \left( \begin {array}{cccccc}
-1&0&0&1&0&0\\ 
0&1&0&0&-1&0\\
0&0&-1&0&0&1\\
0&0&1&0&0&1\\
0&1&0&0&1&0\\
1&0&0&1&0&0
\end {array} \right).
\end{equation}
\item  The set $S^{-1} \Omega  S=\{S^{-1}  A_i  S\}_{i=1}^6$ is  block diagonal matrix where
\begin{footnotesize}
\begin{equation}
\begin{array}{ll}
S^{-1}  A_1  S=
\left(   \begin {array}{ccc}
0&0&0\\
0&0&1\\ 
0&-1&0
\end {array} 
\right) \oplus
\left( \begin {array}{ccc}
0&0&0\\ 
0&0&-1\\
0&1&0
\end {array}
 \right),&
S^{-1}  A_2  S=
\left( \begin {array}{ccc}
0&0&1\\ 
0&0&0\\ 
-1&0&0
\end {array}
\right) \oplus
\left( \begin {array}{ccc}
0&0&-1\\ 
0&0&0\\ 
1&0&0
\end {array}
 \right),\\
S^{-1}  A_3  S=
\left( \begin {array}{ccc}
 0&1&0\\
-1&0&0\\ 
0&0&0
\end {array}
\right) \oplus
\left(  \begin {array}{ccc}
 0&-1&0\\ 
1&0&0\\ 
0&0&0
\end {array}
\right),&
S^{-1}  A_4  S=
\left( \begin {array}{ccc}
0&-1&0\\ 
1&0&0\\
0&0&0
\end {array}
\right) \oplus
\left( \begin {array}{ccc}
0&-1&0\\ 
1&0&0\\ 
0&0&0
\end {array}
 \right),\\
S^{-1}  A_5  S=
\left( \begin {array}{ccc}
0&0&1\\ 
0&0&0\\ 
-1&0&0
\end {array}
\right) \oplus
\left(   \begin {array}{ccc}
0&0&1\\ 
0&0&0\\ 
-1&0&0
\end {array}
\right),&
S^{-1}  A_6 S=
\left( \begin {array}{ccc}
0&0&0\\
0&0&-1\\ 
0&1&0
\end {array} 
\right) \oplus
\left(  \begin {array}{ccc}
0&0&0\\
0&0&-1\\
0&1&0
\end {array}
\right).\\
\end{array}
\end{equation}
\end{footnotesize}
\end{itemize}
\end{example}
\section{Summary}
It is well known that a set of non-defect matrices can be simultaneously diagonalized if and only if the matrices commute. In the case of non-commuting matrices, the best that can be achieved is simultaneous block diagonalization. Both of Algorithm B and Maehara et al \cite{Maehara2011} algorithm are applicable for simultaneous block diagonalization of a set of matrices by a unitary matrix. Algorithm C can be applied for block diagonalization by an invertible matrix when finding a unitary matrix is not possible. In case, block diagonalization of a set of matrices is not possible by both of a unitary or an invertible matrix, then one may utilize block triangularization by Algorithm A.
Algorithms A and B are based on the existence of invariant subspaces however, Algorithm C is based on the existence of a commuting matrix which is not necessarily Hermitian unlike Maehara et al algorithm.
\subsection*{Acknowledgments}
Ahmad Y. Al-Dweik and M. T. Mustafa would like to thank Qatar University for its support and excellent research facilities. R. Ghanam and G. Thompson are grateful to VCU Qatar and Qatar Foundation for their support.\\\\
Conflict of Interest: The authors declare that they have no conflict of interest.
\newpage
\section*{Appendix: Maple procedures}
\lstset{basicstyle=\scriptsize}
\begin{lstlisting}[caption={Step 5 in Algorithm $A$},frame=single]
OrthogonalComplement:= proc(S,N)
uses LinearAlgebra;
local S2;
[seq(cat(e,i),i=1..N)];
Matrix([seq([seq(coeff(S[i],%[j]),j=1..N)],i=1..nops(S))]);
NullSpace(%);
S2:=[seq(add(%[j][i]*cat(e,i),i=1..N),j=1..nops(%))];
end proc:
\end{lstlisting}
\begin{lstlisting}[caption={Step 6 in Algorithm $A$},frame=single]
MatrixProjection:= proc(A,S1,S2)
uses LinearAlgebra;
local S,N,N1,SS,T;
S:=[S1[],S2[]]:
N:=nops(S);
N1:=nops(S1);
[seq(cat(e,i),i=1..N)]:
SS:=Matrix([seq(Vector([seq(coeff(S[i],%[j]),j=1..N)]),i=1..N)]):
[seq(SS^(-1).A[i].SS,i=1..nops(A))]:
T:=map(z->Matrix([ZeroMatrix(N-N1,N1),IdentityMatrix(N-N1)]).z.
Matrix([[ZeroMatrix(N1,N-N1)],[IdentityMatrix(N-N1)]]),%);
end proc:
\end{lstlisting}
\begin{lstlisting}[caption={Steps 8 $\&$ 9 in Algorithm $A$},frame=single]
InvertibleToUnitary:= proc(S,N)
uses LinearAlgebra;
local Q,R,U;
[seq(cat(e,i),i=1..N)];
Matrix([seq(Vector([seq(coeff(S[i],%[j]),j=1..N)]),i=1..nops(S))]);
Q, R := QRDecomposition(%);
U:=Q^+;
end proc:
\end{lstlisting}
\begin{lstlisting}[caption={Steps 2 $\&$ 3 in Algorithm $C$},frame=single]
MatrixCentralizer:= proc(A)
uses LinearAlgebra,ArrayTools;
local n,T,X,kern;
n:= Size(A[1],1);
T:=map(z->KroneckerProduct(IdentityMatrix(n),z)
-KroneckerProduct(z^+,IdentityMatrix(n)),A):
X:=Matrix([seq([T[i]],i=1..nops(T))]);
kern := NullSpace(X);
return [seq(Reshape(kern[k],[n,n]),k=1..nops(kern))];
end proc:
\end{lstlisting}
\begin{lstlisting}[caption={Steps 6 $\&$ 7 in Algorithm $C$},frame=single]
GeneralizedEigenspace := proc(C)
uses LinearAlgebra,ArrayTools;
local n,Evalue,kern,S,U;
n:= Size(C,1);
Evalue:=convert(Eigenvalues(C),set);
kern := [seq(NullSpace((C-Evalue[i]*IdentityMatrix(n))^n),
i=1..nops(Evalue))];
S:=Matrix([seq(op(kern[i]),i=1..nops(kern))]);
return Evalue,kern,S;
end proc:
\end{lstlisting}

\end{document}